\documentclass[11pt]{article}

\usepackage{amssymb,latexsym,tikz,graphicx,eepic,caption}
\usepackage{amsthm}
\usepackage{amsmath,mathrsfs}
\usepackage{amscd,mathtools}
\usepackage{rotating}

\usepackage{amsfonts}

\usepackage[all]{xy}
\usepackage[ngerman,english] {babel}

\usepackage[utf8]{inputenc}

\usetikzlibrary{arrows,shapes}

\usepackage[top=0.7in,bottom=1in,left=1.25in,right=1.25in]{geometry}

\usepackage{verbatim}

\usepackage{setspace}

\usepackage{titlesec}

\newtheorem{theorem}{Theorem}[section]
\newtheorem*{theorem*}{Theorem}

\newtheorem{definition}[theorem]{Definition}

\newtheorem{proposition}[theorem]{Proposition}
\newtheorem{corollary}[theorem]{Corollary}
\newtheorem{example}[theorem]{Example}
\newtheorem{lemma}[theorem]{Lemma}

\numberwithin{equation}{section}

\newtheorem*{corollaryintr}{Corollary}

\usepackage[normalem]{ulem}

\def\Z{\mathbb{Z}}

\def\max#1{\text{max}(#1)}

\newcommand*\xbar[1]{%
	\hbox{%
		\vbox{%
			\hrule height 0.5pt 
			\kern0.5ex
			\hbox{%
				\kern-0.1em
				\ensuremath{#1}%
				\kern-0.1em
			}%
		}%
	}%
}
\def\cSn{\xbar{\mathfrak{S}}_n}

\definecolor{mhcblue}{HTML}{0077CC} 
\definecolor{davidsonred}{HTML}{AC1A2F}

\begin{document}

	\title{Special values of spectral zeta functions and combinatorics: Sturm-Liouville problems
		\footnote{This work is supported by the National Natural Science Foundation of China (Grant 12071371). }}
	\author{Bing Xie\footnote{B. Xie, School of Mathematics and Statistics, Shandong University, Weihai 264209, China.
			Email: {\tt xiebing@sdu.edu.cn}},
		Yigeng Zhao\footnote{Y. G. Zhao, Westlake University, Hangzhou 310024, China.
			Email: {\tt zhaoyigeng@westlake.edu.cn  }},
		and Yongqiang Zhao\footnote{Y. Q. Zhao, Westlake University, Hangzhou 310024, China.
			Email: {\tt  zhaoyongqiang@westlake.edu.cn }}
	}
	\date{\today}

	\maketitle
	
	\begin{abstract}
		In this paper, we apply the combinatorial results on counting permutations with fixed pinnacle and vale sets to evaluate the special values of the spectral zeta functions of Sturm-Liouville differential operators. As applications, we get a combinatorial formula for the special values of spectral zeta functions  and give a new explicit formula for Bernoulli numbers.	
		
		\noindent{\it 2020 MSC numbers}: Primary 05A05,  Secondary 11M06, 34B24
		\vspace{04pt}

		\noindent{\it Keywords}: 
		Spectral zeta functions, Special values, Sturm-Liouville problems, Combinatorics, Permutations, Bernoulli number.
	\end{abstract}

	\section{Introduction}\label{Intro}
	
	Spectral zeta functions for Sturm-Liouville problems are natural generalizations of the classical Riemann zeta function, 
	which are important subjects both in number theory and in spectral theory. As in the Riemann zeta function case, 
	the special values of these spectral zeta functions are of great importance and are usually studied by 
	functional-analysis tools. See, for example, \cite{AGHKLT-21}. In this note, we investigate special values of 
	spectral zeta functions of Sturm-Liouville operators through a combinatorial approach and point out a connection
	between special values and the problem of counting permutations with fixed pinnacles and vales. 
	For the convenience of readers, we start with a detailed review of the Sturm-Liouville problem and its 
	spectral zeta functions.

	Consider the Sturm-Liouville eigenvalue problem with the general separated boundary conditions:
	\begin{align}\label{GeneralS-LProblem}
		\begin{cases}
			Ty:=-y''+qy=\lambda y \;\;{\rm on}\;x\in(a, b),\\
			y(a)\cos\alpha -y'(a)\sin\alpha =0,\\
			y(b)\cos\beta -y'(b)\sin\beta =0,
		\end{cases}
	\end{align}
	where $\alpha,\beta\in[0,\pi)$, $q(x)\in L^1[a, b]$ is a given potential function.
	The real  number (by the self-adjointness of $T$) $\lambda$ is called an eigenvalue if the problem \eqref{GeneralS-LProblem}
	has a nontrivial solution (cf. \cite{Weyl}, \cite[Chapter 3, p.43]{Zettl}).
	
	By the spectral theory of self-adjoint operators, it is known that the eigenvalue problem \eqref{GeneralS-LProblem}
	has only discrete real eigenvalues which are bounded below (cf. \cite[Theorem 13.2]{Weidmann},
	\cite[Theorem 4.3.2]{Zettl}, \cite[p.269]{Walter}).
	Denote $\lambda_k$ the $k$-th eigenvalue of the problem \eqref{GeneralS-LProblem}.  In the case that $q(x)\equiv 0$  and $\alpha=\beta=0$ (i.e., Dirichlet boundary condition), we see easily that its $k$-th eigenvalue is  $\frac{k^2\pi^2}{(b-a)^2}$ with the correspondence eigenfunction $\sin(\frac{k\pi(x-a)}{b-a})$. For the general case, the eigenvalues can be compared with the above special case using Pr\"ufer transformation introduced by  Atkinson and Mingarelli \cite{Atkinson-Mingarelli} and we have the following asymptotic formula (cf. \cite[Lemma 2.1]{Atkinson-Mingarelli}, \cite[Theorem 4.3.1.(7)]{Zettl}):
	\begin{equation}\label{Weyl}
		-\infty<\lambda_1<\cdots<\lambda_k=\frac{k^2\pi^2}{\left( b-a \right)^2}(1+o(1))\to +\infty,
		\; \mathrm{as}\ k\to +\infty.
	\end{equation}

	Thus, suppose for any $k\in\mathbb{N}$, $\lambda_k\not=0$, then we can define the spectral zeta function 
	(cf. \cite[Section 2]{Fucci-Graham-Kirsten}, \cite[(2.25)]{Gesztesy-Kirsten}, \cite[(2.3)]{AGHKLT-19}) associated to 
	the problem \eqref{GeneralS-LProblem} as
	\begin{equation}\label{Zeta-L}
		\zeta_T(s):=\sum_{k=1}^{\infty}\frac{1}{\lambda_k^s}.
	\end{equation}
	Thanks to asymptotic formula \eqref{Weyl}, we know that $\zeta_T(s)$ converges at the half-plane
	$\Re(s)>1/2$.
	Moreover, in the case that  $q$ is smooth in \eqref{GeneralS-LProblem},
	$\zeta_T(s)$ can be meromorphically extended to the whole	complex plane (cf. \cite[Section 4]{Fucci-Graham-Kirsten}, 
	\cite{Gesztesy-Kirsten}).  In general, we do not know whether $\zeta_T(s)$ has a functional equation or not.
	In this paper, we study the problem on the evaluation of special values  $\zeta_T(s)$ at positive integers $s=n$ using the integral of  
	Green functions.
	
	This method naturally involves an interesting combinatorial problem on counting permutations with fixed pinnacle and vale sets.
	To illustrate the main problem, we first consider the following simple example.
	
	\begin{example}\label{Exam}
		Consider the Sturm-Liouville problem
		with Dirichlet boundary conditions
		\begin{equation}\label{special-SL}
			\begin{cases}
				Ty:=-y''=\lambda y \; \;{\rm for}\;x\in(0,\pi),\\
				y(0)=y(\pi)=0.
			\end{cases}
		\end{equation}

		The $k$-th eigenvalue of \eqref{special-SL} is $\lambda_k=k^2$, for any $k\geq 1$
		and up to a factor of  \;$2$, its spectral zeta function is exact the Riemann zeta function
		$$
		\zeta_T(s)=\sum_{k=1}^{+\infty}\frac{1}{k^{2s}}=\zeta(2s).
		$$
		
		Mercer's Theorem (cf. Proposition \ref{Mercer}) tells us that
		\begin{equation}\label{Mercer-Special}
			\zeta_T(n)=\sum_{k=1}^{\infty}\frac{1}{\lambda_k^n}
			=\int_0^\pi\cdots\int_0^\pi G(x_1,x_2)\cdots G(x_n,x_1){\rm d}x_1\cdots{\rm d}x_n,
		\end{equation}
		where $G(s,t)=\frac{1}{\pi}\min\{s,t\}\left(\pi-{\rm max}\{s,t\}\right)$ (cf. \cite[1.6]{AGHKLT-21}) is the Green function of
		\eqref{special-SL} at $\lambda=0$.
		
		For $n=1,2,3$, we may use the symmetry of Green functions to calculate the integral in \eqref{Mercer-Special}. For $n=1$,
		$$
		\zeta_T(1)=\int_0^\pi G(x_1,x_1){\rm d}x_1=\frac{1}{\pi}\int_0^\pi x_1(\pi-x_1){\rm d}x_1=\frac{\pi^{2}}{6}.
		$$
		For $n=2$, since Green function satisfies symmetry, i.e., $G(x_1,x_2)=G(x_2,x_1)$, we divide $[0,\pi]^2$ into $2$ parts: $\{(x_1,x_2): x_1>x_2\}$ and $\{(x_1,x_2): x_1<x_2\}$. Thus
		\begin{align*}
			\zeta_T(2)&=\int_0^\pi\int_0^{x_2} G^2(x_1,x_2){\rm d}x_1 {\rm d}x_2+\int_0^\pi\int^{x_1}_{0} G^2(x_1,x_2){\rm d}x_2{\rm d}x_1\\
			&=\frac{2}{\pi^2}\int_0^\pi\int_0^{x_2} x^2_1(\pi-x_2)^2{\rm d}x_1 {\rm d}x_2 =\frac{\pi^{4}}{90}.
		\end{align*}
		For $n=3$, we divide the integral domain into
		$$
		[0,\pi]^3=\bigcup_{\sigma\in \mathfrak{S}_3}D_{\sigma}:=\bigcup_{\sigma\in \mathfrak{S}_3}\{(x_1,x_2,x_3)| \,x_{\sigma(1)}<x_{\sigma(2)}<x_{\sigma(3)}\},
		$$
		where $\mathfrak{S}_3$ is the permutation group of the set $\{1,2,3\}$. On each $D_{\sigma}$, the integrand is 
		$$
		G_{\sigma}:=\frac{1}{\pi^3}x_{\sigma(2)}(\pi-x_{\sigma(2)})x_{\sigma(1)}^2(\pi-x_{\sigma(1)})^2.
		$$ 
		For $\sigma,\tau\in \mathfrak{S}_3$, it is clear that \[ \int_{D_{\sigma}}G_{\sigma}=\int_{D_{\tau}}G_{\tau}.\]
		Therefore, we may take $\sigma=1$ the identity element and obtain
		$$
		\zeta(6)=3!\int_{D_{1}}G_{1}=\frac{3!}{\pi^3}\int_0^\pi\int_0^{x_3}\int_0^{x_2} x_2(\pi-x_2)x^2_1(\pi-x_1)^2{\rm d}x_1 {\rm d}x_2 {\rm d}x_3
		=\frac{\pi^{6}}{945}.
		$$
	\end{example}
	
	\noindent For $n\geq 4$, however, the symmetry of the integrand is not enough to reduce the whole integral into a multiple of one single piece.  It is obvious that direct calculations would become more and more complicated as $n$ getting larger,  
	which is not practical and pointed out by Ashbaugh et al. \cite{AGHKLT-21} before.  In our approach, the following combinatorial problem on counting permutations  with fixed peaks and valleys naturally arises.	
	
	Let $\sigma=(\sigma(1),\cdots,\sigma(n))$  be a permutation of $[n]=\{1,\cdots, n\}$, with circular conditions:
	$\sigma(0)=\sigma(n)$ and $\sigma(n+1)=\sigma(1)$. The set of all circular $n$-permutations is denoted by $\cSn$.
	We call $\sigma(i)$ is a pinnacle (resp. vale) if  $ \sigma(i-1)<\sigma(i)>\sigma(i+1)$ (resp.  $ \sigma(i-1)>\sigma(i)<\sigma(i+1)$).
	Denote $P(\sigma)$= the set of pinnacles of $\sigma$, and  $V(\sigma)$= the set of vales of $\sigma$. For two disjoint subsets $V$ and $P$ of $[n]$, we let
	\[ C_n(V;P):=\{\sigma\in \cSn|\, V(\sigma)=V; \;  P(\sigma)=P\}.\]
	The pair $(V,P)$ is called $n$-admissible if $C_n(V;P)\neq \emptyset$.
	Let $c_n(V;P):=\sharp C_n(V;P) $ be the number of circular permutations with the pinnacle set $P$ and the vale set $V$.  This class number $c_n(V;P)$ can be explicitly given (see Theorem \ref{orbits.number}).

	Using Mercer's Theorem again, we obtain a formula for  the special values of the spectral zeta function
	of the general Sturm-Liouville problem \eqref{GeneralS-LProblem}. Recall the  Green function (cf. \eqref{Greenfunction}) of \eqref{GeneralS-LProblem}  is
	$$
	G(s,t)=-\frac{1}{W}\psi_-(\min\{s,t\})\psi_+({\rm max}\{s,t\}),\;a\leq s, t\leq b,
	$$
	where the symbols $W, \psi_-$ and  $\psi_+$ can be defined in \eqref{psi+-}-\eqref{Wrons}.
	Note that on each admissible pair $(V,P)$, we could determine the integrand uniquely. Our main result is the following theorem.
	\begin{theorem}\label{Main-Thm}
		For the general Sturm-Liouville problem \eqref{GeneralS-LProblem} and $n\geq 1$, we have
		\begin{align*}
			\zeta_T(n)=\frac{(-1)^nn}{W^n}
			& \sum_{(V,P)}  c_n(V;P)\int_{a}^{b}\int_{a}^{x_n}\cdots\int_{a}^{x_2}
			\prod_{m=1}^{n} \psi^{a_m}_{-}(t_{m}) \psi^{2-a_m}_{+}(t_{m})  {\rm d}x_1\cdots{\rm d}x_n,
		\end{align*}
		where the pair $(V,P)$ runs over all the $n$-admissible pairs and
		$$
		a_m=
		\begin{cases}
			2, & \hbox{if}\; m\in V;\\
			0, & \hbox{if}\; m\in P;\\
			1, & \hbox{otherwise}.
		\end{cases}
		$$
	\end{theorem}
	
	As an application of the above general Theorem, we get a closed combinatorial  formula for Bernoulli numbers at even $B_{2n}$. 
	
	\begin{corollaryintr}[Corollary \ref{Bern}]
		For any integer $n>1$, we have
		\begin{equation*}
			B_{2n}=\frac{(-1)^{n+1}}{2^n(2^{2n}-1)}\sum_{(V,P)}\frac{1}{2^{\sharp P}}\prod_{p \in P\setminus\{n\}} { N_{PV}(p) \choose 2} \prod_{r\in [n]\setminus (P\cup V)}N_{PV}(r)\prod_{m=1}^n\frac{1}{2m+N_{PV}(m+1)},
		\end{equation*}
		where the  pair $(V,P)$ runs over all the $n$-admissible pairs and 
		\begin{equation}\label{NPV}
			N_{PV}(\alpha)=\sharp \{ v \in V: v<\alpha\}-\sharp\{ p \in P: p<\alpha\}\; {\rm for\; any}\; 1\leq\alpha\leq n.
		\end{equation}
	\end{corollaryintr}

	The rest of this paper is organized as follows.	In Section \ref{Mercer Thm}, we quickly recall Green functions and Mercer's Theorem.
	We give the explicit formula for the class number $c_n(V;P)$ in Section \ref{Combinatorics} and prove our main theorem in Section \ref{proof}.
	In Section \ref{example}, we give a formula of Riemann zeta function $\zeta(2n)$ and Bernoulli numbers as applications.

	\section{Green functions and Mercer's Theorem}\label{Mercer Thm}

	We can define the domain of the operator $T$ in \eqref{GeneralS-LProblem} as
	$$
	\mathcal{D}(T):=\Big\{u,\,pu'\in AC[a,b]:\,T u\in L^2[a,b]\,{\rm and}\; u\;{\rm satisfies\, the\, boundary\, condition\, in}\, 
	\eqref{GeneralS-LProblem} \Big\},
	$$
	where $AC[a,b]$ represents all the absolutely continuous functions on $[a,b]$
	and
	$$
	L^2[a,b]:=\bigg\{u\;{\rm is\;measurable}:\;\int_a^b |u|^2 {\rm d}x<\infty\bigg\}.
	$$
	Then $T$ is a self-adjoint operator with domain $\mathcal{D}(T)$.	
	
	For any $\lambda_0\in\mathbb{C}$ which is not an eigenvalue of $T$, then the inverse operator of $T-\lambda_0 \mathrm{Id}$,
	$(T-\lambda_0 \mathrm{Id})^{-1}$, exists and is a bounded linear operator on $L^2[a,b]$, where $ \mathrm{Id}$ is the identity operator on $L^2[a,b]$. 
	Then the Green function $G(s,t)$ of the operator $T_1$ at $\lambda_0$ is defined as follows, for any $f\in L^2[a,b]$,
	\begin{equation}\label{T-Green}
		(T-\lambda_0\mathrm{Id})^{-1}f(s)=\int_a^bG(s,t)f(t){\rm d}t,\;s\in[a,b],
	\end{equation}
	and satisfies the boundary condition, i.e., for  any fixed $t\in[a,b]$,
	\begin{equation}\label{eq:bdy-Green}
		\cos\alpha G(a,t)-\sin\alpha \left(\frac{ \partial G }{ \partial s}\right)(a,t)=0
		=\cos\beta G(b,t)-\sin\beta \left( \frac{ \partial G }{ \partial s}\right)(b,t).
	\end{equation}

	The definition \eqref{T-Green}  is equivalent to that, for any fixed $t\in[a,b]$,
	\begin{equation}\label{T-Green-2}
		(T-\lambda_0 \mathrm{Id})G(s,t)=\delta_t(s)
	\end{equation}
	where $\delta_t(s)=\delta(s-t)(=\delta(t-s))$ is the Dirac impulse at $t$, and it is defined via
	$$
	\int_a^b \delta_t(s)f(s){\rm d}s=f(t),
	$$
	for any $t\in[a,b]$ and continuous function $f(s)$ on $[a,b]$.
	
	Recall the Mercer's Theorem for Sturm-Liouville eigenvalue problems \eqref{GeneralS-LProblem} is given as follows.
	
	\begin{proposition} \cite[\S 3.5.4]{Courant-Hilbert}(Mercer's Theorem)\label{Mercer}
		
		For $\lambda_0\not=\lambda_k$, $k\ge 1$, we denote $G(x,t),\,x,\,t\in [a,b]$, the Green function
		of Problem \eqref{GeneralS-LProblem}  at $\lambda_0$. Then for any integer $n\ge 2$,
		\begin{equation}\label{Mercer1}
			\int_a^b\cdots\int_a^bG(x_1,x_2)\cdots G(x_n,x_1){\rm d}x_1\cdots{\rm d}x_n
			=\sum_{k=1}^{\infty}\frac{1}{(\lambda_k-\lambda_0)^n}.
		\end{equation}
	\end{proposition}
	
	Indeed, In \cite[\S 3.5.4]{Courant-Hilbert}, Courant and Hilbert have given the proof of Mercer's theorem for 
	linear integral operator. We note that a solution of the Sturm-Liouville problem \eqref{GeneralS-LProblem} can be written in the form
	\begin{equation}
		y(x)=(\lambda-\lambda_0)\int_{a}^{b}G(x,t)y(t){\rm d}t,
	\end{equation}
	which is an integral equation with kernel $G(x,t)$.  Therefore Proposition \ref{Mercer} can be deduced from  \cite[\S 3.5.4]{Courant-Hilbert}. Moreover, as  \cite[(3.6)]{AGHKLT-21} shown, Proposition \ref{Mercer} also holds for trace class operators.

	By translating the potential function $q$, we may suppose that $\lambda=0$ is not an eigenvalue of $T$.
	In fact, if there exists an eigenvalue $\lambda_k=0$, we can consider replacing problem \eqref{GeneralS-LProblem} with
	a translated problem 
	\begin{align}\label{S-LProblem-trans}
		\begin{cases}
			-y''+\left(q-\frac{1}{2}\lambda_{k+1}\right)y=\lambda y \;\;{\rm on}\;x\in(a, b),\\
			y(a)\cos\alpha -y'(a)\sin\alpha =0,\\
			y(b)\cos\beta -y'(b)\sin\beta =0.
		\end{cases}
	\end{align}
	Clearly, the $k$-th and $(k+1)$-th eigenvalues of \eqref{S-LProblem-trans} are $-\frac{1}{2}\lambda_{k+1}$ and $\frac{1}{2}\lambda_{k+1}$,
	respectively.
	
	Hence, in the following and in Proposition \ref{Mercer}, without loss of generality, we can always suppose $\lambda_0=0$.
	Therefore, we get an integral representation of spectral zeta function $\zeta_T(s)$ at positive integers.
	More precisely,	for any $n\ge 1$, we have
	\begin{equation}\label{Spec-Mercer}
		\zeta_T(n)=\sum_{k=1}^{\infty}\frac{1}{\lambda_k^n}
		=\int_a^b\cdots\int_a^bG(x_1,x_2)\cdots G(x_n,x_1)){\rm d}x_1\cdots{\rm d}x_n.
	\end{equation}
	
	Following the method in \cite[(10), p.252]{Walter}, 
	in order to find the Green function, we put $\psi_-$ and $\psi_+$ be two solutions of differential
	equation $-y''+qy=0$, (i.e., \eqref{GeneralS-LProblem} with $\lambda=0$) with Cauchy conditions
	\begin{equation}\label{psi+-}
		\psi_-(a)=\sin\alpha,\; \psi_-'(a)=\cos\alpha, \; \mathrm{and}\; \psi_+(b)=\sin\beta,\;\psi'_+(b)=\cos\beta,
	\end{equation}
	respectively. Since we suppose that $\lambda=0$ is not an eigenvalue of problem \eqref{GeneralS-LProblem},
	$\psi_-$ and $\psi_+$ are real value linear independence solutions. This fact leads to	
	\begin{align}
		W& :=W[\psi_-, \psi_+]:=\psi_-(x)\psi_+'(x)-\psi_+(x)\psi_-'(x)\nonumber\\
		&=\sin\alpha\psi_+'(a)-\psi_+(a)\cos\alpha\neq 0,  \label{Wrons}
	\end{align}
	where the constant $W=W[\psi_-, \psi_+]$ is the Wronskian of the solutions $\psi_-$ and  $\psi_+$. 
	
	The definition of Green function  \eqref{T-Green-2} at the point $\lambda_0=0$ is equivalent to that, for any fixed $t\in[a,b]$,
	$G(s,t)$ is continuous and satisfies 
	\begin{equation}\label{Green-temp1}
		-\frac{ \partial^2 G }{ \partial s^2}(s,t)+q(s)G(s,t)=0,\, \;\mathrm{for}\; s\not=t,
	\end{equation}
	and 
	\begin{equation}\label{Green-temp2}
		\frac{ \partial G }{ \partial s}(t-,t)-\frac{ \partial G }{ \partial s}(t+,t)=1.
	\end{equation}
	with the boundary condition \eqref{eq:bdy-Green}.
	By the definition of  $\psi_-$ and $\psi_+$, we see the Green function should be 
	$$
	G(s,t)=\begin{cases}
		C_1(t)\psi_-(s), & a\leq s\leq t\leq b, \\
		C_2(t)\psi_+(s), & a\leq t\leq s\leq b,
	\end{cases}
	$$
	where $C_1(t)$ and $C_2(t)$ are functions that can be determined in the following. Clearly, $G(s,t)$ satisfies \eqref{Green-temp1} and \eqref{eq:bdy-Green},
	by the boundary conditions \eqref{psi+-}. 
	
	Due to the continuity of $G(s,t)$ at $s=t$, we have $C_1(t)\psi_-(t)=C_2(t)\psi_+(t)$, 
	hence $C_1(t)=C\psi_+(t)$ and $C_2(t)=C\psi_-(t)$, for some constant $C$. Using \eqref{Green-temp2}, we  finally see that $C=-1/W$. In summary,  the Green function of problem \eqref{GeneralS-LProblem} at the point $\lambda_0=0$ is given by
	\begin{align}
		G(s,t)&=-\frac{1}{W} 
		\begin{cases}
			\psi_-(s)\psi_+(t), & a\leq s\leq t\leq b, \\
			\psi_-(t)\psi_+(s), & a\leq t\leq s\leq b,
		\end{cases}
		\nonumber \\
		&=-\frac{1}{W}\psi_-(\min\{s,t\})\psi_+({\rm max}\{s,t\}),\;a\leq s, t\leq b.  \label{Greenfunction}
	\end{align}

	To calculate the special value $\zeta_T(n)$, we have to divide the range of the integral $[a,b]^n $ to  determine 
	the integrand $G(x_1,x_2)\cdots G(x_n,x_1)$ in \eqref{Spec-Mercer}. Thanks to the explicit formula \eqref{Greenfunction}, 
	this leads to the necessity to compare the values of $x_i$ with its adjacencies $x_{i+1}$ and $x_{i-1}$ for $1\leq i\leq n$, 
	where we take $x_0=x_n$ and $x_{n+1}=x_1$ for our convenience. More precisely, given $(x_1,\cdots,x_n)\in [a,b]^n$, 
	we see the term $\psi_{-}^{a_i}(x_i)\psi_{+}^{2-a_i}(x_i) $ involved $x_i$ , for $1\leq i\leq n$,  appears in the integrand, where
	\begin{equation}\label{symbol}
		a_i=\begin{cases*}
			2, & if $x_{i-1}>x_i<x_{i+1}$;\\
			0, &  if  $x_{i-1}<x_i>x_{i+1}$;\\
			1,&  otherwise.
		\end{cases*}
	\end{equation}
	Using the symbol \eqref{symbol}, we have
	\begin{align}\label{Green-count}
		(-W)^{n}G(x_1,x_2)\cdots G(x_n,x_1) &=(-W)^{n}\prod_{i=1}^{n} G(x_i,x_{i+1}) \nonumber \\
		=& \prod_{i=1}^{n} \psi_{-}(\min\{x_{i}, x_{i+1}\})\psi_{+}(\min\{x_{i}, x_{i+1}\}) \nonumber\\
		=&\prod_{i=1}^{n} \psi^{a_i}_{-}(x_i)  \psi^{2-a_i}_{+}(x_i).
	\end{align}
	This equality plays a key role in the calculation of the spectral zeta function \eqref{Spec-Mercer}.

	The comparison of the values of $x_i$ may also be reformulated as following: given any $(x_1,x_2\cdots,x_n)\in [a,b]^n$, we may assume 
	$x_i\not= x_j$, for any $i\not=j$. Then there exists a unique $\sigma \in \mathfrak{S}_n$, such that 
	$$
	x_{\sigma(1)}< x_{\sigma(2)}< \cdots < x_{\sigma(n)}.
	$$
	Set $t_i:=x_{\sigma(i)}$ and  $\tau:=\sigma^{-1} \in \mathfrak{S}_n$. Then
	\begin{equation}\label{t_i}
		t_1<t_2<\cdots<t_n,\;x_i=x_{(\sigma^{-1}\sigma)(i)}=t_{\tau(i)},
	\end{equation}
	and hence
	\begin{equation}\label{Green-change}
		G(x_1,x_2)\cdots G(x_n,x_1)=G(t_{\tau(1)},t_{\tau(2)})\cdots G(t_{\tau(n)},t_{\tau(1)}).
	\end{equation}
	We can find that permutations $\tau \in \mathfrak{S}_n$ satisfies $x_i=t_{\tau(i)}< t_{\tau(j)}=x_j$  
	if and only if $\tau(i)<\tau(j)$ for any $1\leq i,j \leq n$. 
	Therefore, comparing the values of $x_i=t_{\tau(i)}$ is equivalent to comparing the values of $\tau(i)$ and 
	this arises naturally a combinatorial question on permutations.

	\section{Combinatorics}\label{Combinatorics}
	We have seen in the last section, the integrand $G(x_1,x_2)\cdots G(x_n,x_1)$ in \eqref{Spec-Mercer} 
	is determined by the numbers $a_i$'s in \eqref{symbol}.  This computation is reduced to find pinnacles and vales in a permutation, 
	which is a classical topic in combinatorial mathematics and has been well-studied in the literature, 
	for example, see \cite{BBS,DNPT,BCMY,DLHHIN,Entringer} and the references there.
	Nevertheless, most of the results on counting peaks and valleys are for usual permutations. 
	As the special form of our integrand above, we require that for circular permutations. 
	
	\begin{definition}
		\begin{itemize}
			\item[(i)] A sequence $\sigma=(\sigma(1),\cdots,\sigma(n))$  is called a circular $n$-permutation if it is a permutation of $[n]=\{1,\cdots,n\}$ and circular: $\sigma(0)=\sigma(n), \sigma(n+1)=\sigma(1)$.  The set of all circular $n$-permutations is denoted by $\cSn$. As a quotient of the group of $n$-permutations $\mathfrak{S}_n$, $\cSn=\mathfrak{S}_n/(\Z/n\Z)$, where $\Z/n\Z$ acts on $\mathfrak{S}_n$ by the natural rotation.
			\item[(ii)] For  $\sigma=(\sigma(1) ,\cdots,\sigma(n)) \in \cSn$, we call $\sigma(i)$ is a pinnacle (resp. vale) if  $ \sigma(i-1)<\sigma(i)>\sigma(i+1)$ (resp.  $ \sigma(i-1)>\sigma(i)<\sigma(i+1)$). We denote $P(\sigma)$= the set of pinnacles of $\sigma$, and  $V(\sigma)$= the set of vales of $\sigma$.
		\end{itemize}
	\end{definition}
	\begin{example}\label{1n}
		For a circular $n$-permutation, $1$ is always a vale and $n$ is always a pinnacle. If $n\geq 3$, $2$ cannot be a pinnacle  and $n-1$ cannot be a vale. For example $\sigma=(2314576) \in \xbar{\mathfrak{S}}_7$, then $P(\sigma)=\{3,7\}$ and $V(\sigma)=\{1,2\}$.
	\end{example}

	\begin{lemma}
		For a circular $n$-permuation, the numbers of pinnacles and vales are the same.
	\end{lemma}
\begin{proof}
	It suffices to notice that any two adjacent pinnacles have a vale between them.
\end{proof}
In this section, we are interested in counting the number of $n$-permutations with fixed pinnacles and vales.
\begin{definition}\label{class-number}
	\begin{itemize}
		\item[(i)]	For $k\in \mathbb{N}$ and  $n\geq 2k$, $V$ and $P$ are  two disjoint subsets with $k$ elements of positive integers.
		We define
		\[ C_n(V;P)=\{\sigma\in \cSn| V(\sigma)=V; \;  P(\sigma)=P\}\]
		\[ c_n(V;P)=\sharp C_n(V;P). \]
		We take  $c_n(V;P) = 0$ if $C_n(V;P)=\emptyset$ for our convenience.
		\item[(ii)] Two disjoint subsets $(V,P)$ of positive integers are called $n$-admissible if  $C_n(V,P)\neq \emptyset$.  We say $(V,P)$ is admissible if it is $n$-admissible for some $n\in \mathbb{N}$.
	\end{itemize}
	
\end{definition}

As $1$ is always a vale and $n$ is always a pinnacle for an $n$-permutation, to simplify our notations, we reformulate the definition of admissible pairs as follows.
\begin{definition}
Let $Q$ and $W$ be two sets of $\mathbb Z_{\geq 2}$.
The pair $(Q,W)$ is called $n$-admissible if $n > \max Q$ and there exists $\sigma \in \cSn$ such that $P(\sigma) = Q \cup \{n\}, V(\sigma) = W \cup \{1\}$. We say $(Q,W)$ is admissible if $(Q,W)$ is $n$-admissible for some $n$.
\end{definition}
By inserting $n+1$ adjacent to $n$ in the permutation $\sigma$, we see that $(Q,W)$ is $n$-admissible implies that $(Q,W)$ is $n+1$-admissible. By induction, it is $m$-admissible for any $m\geqslant n$.

Note that the necessary conditions of a pair $(Q,W)$  to be admissible are $|Q|=|W|$ and $Q\cap W=\emptyset$. For simplification, we always assume these conditions. Given $(Q,W)$, we may write $Q=\{q_1 < \cdots < q_k\}, W= \{w_1 < \cdots < w_k\}$.

\begin{proposition}
The pair $(Q,W)$ is admissible  if and only if $q_j > w_j, \text{ for any } j \in [k]$.	
\end{proposition}	
\begin{proof}
If $(Q,W)$ is admissible, then there exists $\sigma \in \cSn$ such that $P(\sigma) = Q \cup \{n\}$ and $V(\sigma) = W \cup \{1\}$.   We note that,

For $w_1\in W$, if $q_1$ is one of the adjacent pinnacles of $w_1$, then it is clear $q_1>w_1$. Otherwise, $q_1$ has two adjacent vales in $(W\setminus\{w_1\})\cup \{1\}$, say $w'\in W$ is one of that, then we have $q_1>w'>w_1$.  Now for $w_2\in W$,  there are three possibilities:
\begin{itemize}
	\item[(i)]  if $w_2$ is one of the adjacent vales of $q_2$, it is clear that $q_2>w_2$;
	\item[(ii)] if  $w_1$ is one of the adjacent vales of $q_2$, then we have $q_2>q_1>w'\geq w_2$. 
	\item[(iii)] if neither $w_1$ nor $w_2$ is one of the adjacent vales of $q_2$, then it follows there exists $w''\in (W\setminus\{w_1,w_2\})\cup \{1\} $ such that $q_2>w''$. Hence $q_2>w''>w_2$.
\end{itemize}
We continue this process and  obtain the assertion inductively.

Conversely, using $Q$ and $W$, we can define a word $\sigma$ as follows:
\begin{equation*}
	\sigma(i)= \begin{cases}
		1 & \mathrm{if}\; i=1 \\
		q_j & \mathrm{if}\; i=2j \; \mathrm{for}\; 1\leq j\leq k\\
		w_j &\mathrm{if}\; i=2j+1 \; \mathrm{for}\; 1\leq j\leq k\\
		q_k+1 & \mathrm{if}\;  i=2k+2
	\end{cases}
\end{equation*}
Now we want to insert elements into the word $\sigma$ to produce a permutation.

For $m=\max Q +1 = q_k +1$. We denote $T=[m]\setminus (Q \cup W \cup \{1,m\})$, and insert  the elements of $T$ between $\sigma(2k+2)=q_k+1$ and $\sigma(1)=1$ in the descending order. In this way, we obtain an element  $\tilde{\sigma} \in \xbar{\mathfrak{S}}_m$ with $P(\tilde{\sigma}) = Q \cup \{m\}, V(\tilde{\sigma}) = W \cup \{1\}$. Therefore, $(Q,W)$ is $\max Q +1$-admissible.
\end{proof}

\begin{corollary}
If $(Q,W)$ is admissible, then $(Q\setminus\{q_k,\cdots, q_{k-i}\}, W\setminus\{w_k,\cdots,w_{k-i}\})$ is also admissible, for any $i\in [k]$.
\end{corollary}

A variant of the class number $c_n(V;P) $ for non-circular permutations has been studied in \cite{DLHHIN}. Its explicitly formula can be  used in our case by a simple modification. In fact, the Foata-Strehl \cite{FS,Braenden} action also works well for circular permutations,  and the method in \cite{DLHHIN} can be also adapted to our setting. We quote the following theorem of Diaz-Lopez et al. directly without calculating $c_n(V;P) $  explicitly in $\cSn$ to keep this section concise and coherent.
\begin{theorem}({\cite[Theorem 4.4]{DLHHIN}})\label{orbits.number}
Let $(V,P)$ be an $n$-admissible pair.	We have the class number \[ c_n(V;P)=2^{n-k-1}\prod_{p \in P\setminus\{n\}} { N_{PV}(p) \choose 2} \prod_{r\in [n]\setminus (P\cup V)}N_{PV}(r), \]
where $k=\sharp P=\sharp V$, and $N_{PV}(\alpha)=\sharp \{ v \in V: v<\alpha\}-\sharp\{ p \in P: p<\alpha\}$ for any $1\leq\alpha\leq n$.
\end{theorem}
\begin{proof}
The permutations considered in \cite{DLHHIN} are not circular, but we may write our circular permutation $\sigma=(\sigma(1),\cdots,\sigma(n)) =(n \tilde{\sigma})$ with $\tilde{\sigma}$ is a usual $(n-1)$-permutation, then we apply the main result there for the pinnacles set $P\setminus\{n\}$ in the first product.
\end{proof}
\begin{example}
For $n=6$, we have the following table of the class number $c_6(V;P)$ for admissible pairs $(V,P)$:
\begin{center}
	\begin{tabular}{|c|c|c|c|c|c|c|c|c|c|}
		\hline
		V &  $\{1\}$ &$\{1,2\}$& $\{1,2\}$& $\{1,2\}$&$\{1,3\} $&$\{1,3\}$&$\{1,4\}$&$\{1,2,3\}$&$\{1,2,4\}$ \\
		\hline
		P&  $\{6\} $&$\{3,6\}$& $\{4,6\}$&$\{5,6\} $&$\{4,6\}$&$\{5,6\}$&$\{5,6\}$&$\{4,5,6\}$&$\{3,5,6\} $\\
		\hline
		$c_6(V;P)$ &$16$&$8$&$16$&$32$&$8$&$16$&$8$&$12$&$4$\\
		\hline
	\end{tabular}
\end{center}

\end{example}
\section{The proof of main theorem}\label{proof}
Now, we  give {the proof of Theorem \ref{Main-Thm}}.  In fact, with the representation of the combination number $c_n(V;P)$,
the proof implicit in Example \ref{Exam}. 

Recall some symbols in case $n=3$ of Example \ref{Exam} and \eqref{t_i}. 
For any $n\geq 1$, we divide the integral domain into
\begin{equation}\label{ab-D}
[a,b]^n=\bigcup_{\sigma\in \mathfrak{S}_n}D_{\sigma}:=\bigcup_{\sigma\in \mathfrak{S}_n}\{(x_1,x_2,\cdots,x_n)| 
\,x_{\sigma(1)}<x_{\sigma(2)}<\cdots <x_{\sigma(n)}\},\;a.e..
\end{equation}
Set $t_i:=x_{\sigma(i)}$, $\tau:=\sigma^{-1}\in\mathfrak{S}_n$, and hence $x_i=t_{\sigma^{-1}(i)}=t_{\tau(i)}$ and 
$t_1<t_2<\cdots<t_n$.
Then  by \eqref{symbol}-\eqref{Green-change}, we have  that  on each 
\begin{align*}
D_{\sigma}&=\{(x_1,x_2,\cdots,x_n)| \,x_{\sigma(1)}<x_{\sigma(2)}<\cdots <x_{\sigma(n)}\}\\
&=\{(t_{\tau(1)},t_{\tau(2)},\cdots,t_{\tau(n)})| \,t_1<t_2<\cdots<t_n\}
\end{align*}
the integrand is 
\begin{align}\label{Green-n}
G_{\sigma}:= & G(x_1,x_2)\cdots G(x_n,x_1) \nonumber\\
= & G(t_{\tau(1)},t_{\tau(2)})\cdots G(t_{\tau(n)},t_{\tau(1)})  \nonumber\\
=&\frac{(-1)^n}{W^n}\psi^{a_{\tau(1)}}_{-}(t_{\tau(1)})  \psi^{a_{\tau(2)}}_{-}(t_{\tau(2)})\cdots \psi^{a_{\tau(n)}}_{-}(t_{\tau(n)}) \nonumber\\
\;\;\;\;\;\; & \cdot \psi^{2-a_{\tau(1)}}_{+}(t_{\tau(1)})\psi^{2-a_{\tau(2)}}_{+}(t_{\tau(2)})\cdots \psi^{2-a_{\tau(n)}}_{+}(t_{\tau(n)})\nonumber\\
=&\frac{(-1)^n}{W^n}\psi^{a_1}_{-}(t_{1})  \psi^{a_2}_{-}(t_{2})\cdots \psi^{a_{n}}_{-}(t_{n})
\psi^{2-a_1}_{+}(t_{1})\psi^{2-a_2}_{+}(t_{2})\cdots \psi^{2-a_{n}}_{+}(t_{n}) \nonumber\\
=&\frac{(-1)^n}{W^n}\prod_{m=1}^{n} \psi^{a_m}_{-}(t_{m}) \psi^{2-a_m}_{+}(t_{m}). 
\end{align}
By \eqref{symbol} and $t_{\tau(i)}< t_{\tau(j)}$ if and only if $\tau(i)<\tau(j)$ for any $1\leq i,j \leq n$,
we have
\begin{equation}\label{symbol-temp}
a_{\tau(i)}=\begin{cases*}
	2, & if $\tau(i-1)>\tau(i)<\tau(i+1)$;\\
	0, &  if $\tau(i-1)<\tau(i)>\tau(i+1)$;\\
	1,&  otherwise.
\end{cases*}
\end{equation}

Suppose $\tau\in C_n(V,P)$, where the pair $(V,P)$ is an $n$-admissible pair. 
By \eqref{symbol-temp} we know that if $\tau(i)\in V$ (resp.  $\tau(i)\in P$), 
then the term $\psi_{-}^2(t_{\tau(i)})$ (resp. $\psi_{+}^2(t_{\tau(i)})$) 
appears in the integrand. Otherwise the term $\psi_{-}(t_{\tau(i)})\psi_{+}(t_{\tau(i)})$ appears.
Hence set $m=\tau(i)$, then \eqref{symbol-temp} can be rewritten as
\begin{equation}\label{symbol-am}
a_m=
\begin{cases}
	2, & \hbox{if}\; m\in V;\\
	0, & \hbox{if}\; m\in P;\\
	1, & \hbox{otherwise}.
\end{cases}
\end{equation}
Moreover, using \eqref{Green-n} and \eqref{symbol-am}, we can calculate the integral of the Green function $G_{\sigma}$
on $D_{\sigma}$,

\begin{align}\label{Int-G}
&\int_{D_{\sigma}}G_{\sigma}{\rm d}x_{1}\cdots{\rm d}x_{n}=
\frac{(-1)^n}{W^n}\int_a^b\int_a^{t_n}\cdots\int_a^{t_2}
\prod_{m=1}^{n} \psi^{a_m}_{-}(t_{m}) \psi^{2-a_m}_{+}(t_{m}){\rm d}t_{1}\cdots{\rm d}t_{n}.
\end{align}
Since $t_1<t_i$, for any $i>1$ and $t_n>t_j$ for any $j<n$, we can get $t_1\in V$ and $t_n\in P$. Hence $a_1=2$ and $a_n=0$ 
in \eqref{Int-G}.

From \eqref{symbol-am}, we can see that the formula \eqref{Int-G} is determined by the admissible pair $(V,P)$.
Hence, for any $n$-admissible pair $(V,P)$ and any two permutations $\tilde{\tau}=\tilde{\sigma}^{-1},\tau=\sigma^{-1} \in C_n(V,P)$, 
it is clear that 
\begin{equation}\label{int=}
\int_{D_{\sigma}}G_{\sigma}{\rm d}x_{1}\cdots{\rm d}x_{n}
=\int_{D_{\tilde{\sigma}}}G_{\tilde{\sigma}}{\rm d}x_{1}\cdots{\rm d}x_{n}.
\end{equation}

With these preparations, we can calculate 	$\zeta_T(n)$.	
By the integral expression \eqref{Spec-Mercer} of $\zeta_T(n)$, \eqref{ab-D}  and  $\cSn=\mathfrak{S}_n/(\Z/n\Z)$,
we have 

\begin{align*}
\zeta_T(n) = & \sum_{\sigma\in\mathfrak{S}_n} \int_{D_{\sigma}}G_{\sigma}{\rm d}x_{1}\cdots{\rm d}x_{n} =n\sum_{\sigma\in\cSn} \int_{D_{\sigma}}G_{\sigma}{\rm d}x_{1}\cdots{\rm d}x_{n}\\
=&n\sum_{(V,P)}\sum_{\sigma\in C_n(V,P)} \int_{D_{\sigma}}G_{\sigma}{\rm d}x_{1}\cdots{\rm d}x_{n}\\
=&n\sum_{(V,P)} c_n(V,P) \int_{D_{\sigma}}G_{\sigma}{\rm d}x_{1}\cdots{\rm d}x_{n}.
\end{align*}
The last two equalities are due to $\cSn=\bigcup_{(V,P)} C_n(V,P)$ and \eqref{int=},
where the pair $(V,P)$ runs over all the $n$-admissible pairs.

Then, by substituting \eqref{Int-G} into the last equality, we can obtain
\begin{align*}
\zeta_T(n) =&\frac{(-1)^nn}{W^n}  \sum_{(V,P)}  c_n(V;P)
\int_a^b\int_a^{t_n}\cdots\int_a^{t_2} 
\prod_{m=1}^{n} \psi^{a_m}_{-}(t_{m}) \psi^{2-a_m}_{+}(t_{m}) {\rm d}t_{1}\cdots{\rm d}t_{n},
\end{align*}
where the pair $(V,P)$ runs over all the $n$-admissible pairs and $a_m$ is in \eqref{symbol-am}.

The proof of Theorem \ref{Main-Thm} is finished. \qed

In Theorem \ref{Main-Thm}, by using the two independent solutions of Sturm-Liouville differential equation and combining the combination numbers
generated by the pinnacles and vales of permutations, an explicit expression of the spectral zeta function $\zeta_T(n)$ is obtained.

\section{An example }\label{example}

In 1735, L. Euler first found the values $\zeta(2n)$ for all even numbers by using Bernoulli numbers.
In this section, we apply our main Theorem  \ref{Main-Thm}  to an explicit example. 
As a byproduct, we give another way to evaluate $\zeta(2n)$, $n\ge 1$ and 
obtain a new closed combinatorial formula for Bernoulli numbers.

\begin{example}\label{Example}
Consider the following special case of Sturm-Liouville problem \eqref{GeneralS-LProblem}:
\begin{equation}\label{SpecialS-LProblem}
	\begin{cases}
		Tu:=-u''=\lambda u \;\; {\rm for}\; x\in [0,1],\\
		u(0)=u'(1)=0.
	\end{cases}	
\end{equation}
Then the $k$-th eigenvalue of problem \eqref{SpecialS-LProblem} is $\lambda_k=\left(k-1/2\right)^2\pi^2$,
and the corresponding eigenfunction is
$$
\sin(\sqrt{\lambda_k} x)=\sin\left((k-1/2)\pi x\right),\;x\in [0,1].
$$
Then, the spectral zeta function of \eqref{SpecialS-LProblem} is
$$
\zeta_T(s)=\sum_{k=1}^{\infty}\frac{1}{(k-1/2)^{2s}\pi^{2s}}.
$$

The $\psi_-$ and $\psi_+$ of problem \eqref{SpecialS-LProblem} are $\psi_-(x)=x$, $\psi_+(x)=1$.
Hence the Wronskian of these two solutions $\psi_-$ and  $\psi_+$ is
$$
W=W[\psi_-, \psi_+]=\psi_-(x)\psi_+'(x))-\psi_+(x)\psi_-'(x)=-1.
$$

\noindent Using these facts, by \eqref{Greenfunction}, we can get the Green function of problem \eqref{SpecialS-LProblem} 
at $\lambda_0=0$ is 
(\cite[\S 5.15.1]{Courant-Hilbert}),
\begin{equation}\label{Green-temp}
	G(s,t)=\min\{s, t\}=
	\begin{cases}
		s, & 0\leq s\leq t\leq 1, \\
		t, & 0\leq t\leq s\leq 1.
	\end{cases}
\end{equation}

Then by Theorem \ref{Main-Thm}, we can get
\begin{align*}
	\zeta_{T}(n)=&n\sum_{(V,P)}c_n(V,P)\int_{0}^{1}\int_{0}^{x_n}\cdots\int_{0}^{x_2}x_1^2x_2^{a_2}\cdots x_{n-1}^{a_{n-1}}{\rm d}x_1\cdots{\rm d}x_n\\
	=&n\sum_{(V,P)}c_n(V,P)\frac{1}{3}\frac{1}{4+a_2}\cdots\frac{1}{m+2+a_2+\cdots +a_m}
	\cdots\frac{1}{n+2+a_2+\cdots +a_{n-1}+0}\\
	=&n\sum_{(V,P)}c_n(V,P)\prod_{m=1}^n\frac{1}{m+a_1+\cdots+a_m}\\
	=&n\sum_{(V,P)}c_n(V,P)\prod_{m=1}^n\frac{1}{2m+N_{PV}(m+1)},
\end{align*}
where $\sum_{(V,P)}$ means the summation of all the $n$-admissible pairs $(V,P)$,   and
$$
a_m=
\begin{cases}
	2, & \hbox{if}\; m\in V;\\
	0, & \hbox{if}\; m\in P;\\
	1, & \hbox{otherwise}.
\end{cases}
$$
Indeed, the last equality follows from the definition of $N_{PV}(\alpha)$, 
\begin{equation}\label{NPV-temp}
	N_{PV}(\alpha)=\sharp \{ v \in V: v<\alpha\}-\sharp\{ p \in P: p<\alpha\}\; {\rm for\; any}\; 1\leq\alpha\leq n.
\end{equation}
Then we can obtain $N_{PV}(m+1)=a_1+\cdots+a_m-m$.

Plugging in the class number formula $c_n(V;P)$,  we get
\begin{align*}\label{zeta(2n)}
	\zeta_{T}(n)=2^{n-1}n\sum_{(V,P)}\frac{1}{2^{\sharp P}}\prod_{p \in P\setminus\{n\}} { N_{PV}(p) \choose 2} \prod_{r\in [n]\setminus (P\cup V)}N_{PV}(r)\prod_{m=1}^n\frac{1}{2m+N_{PV}(m+1)}.
\end{align*}
\end{example}
\vskip 4mm

\noindent We note that 
\begin{align*}
\zeta(s)&=\sum_{k=1}^{\infty}\frac{1}{k^{s}}=\sum_{k=1}^{\infty}\frac{1}{(2k-1)^{s}}+
\sum_{k=1}^{\infty}\frac{1}{(2k)^{s}}
=\frac{1}{2^{s}}\sum_{k=1}^{\infty}\frac{1}{(k-1/2)^{s}}+\frac{1}{2^{s}}\sum_{k=1}^{\infty}\frac{1}{k^{s}} \nonumber  \\
&=\frac{1}{2^{s}}\sum_{k=1}^{\infty}\frac{\pi^s}{(k-1/2)^{s}\pi^s}+\frac{1}{2^{s}}\zeta(s)
=\frac{\pi^s}{2^{s}}\zeta_T(s/2)+\frac{1}{2^{s}}\zeta(s).
\end{align*}
Hence as $\Re(s)>1$, we have
\begin{equation}\label{Zeta-temp}
\zeta(s)=\frac{\pi^s}{2^{s}-1}\zeta_T(s/2).
\end{equation}

\noindent Therefore,  
$$\zeta(2n)=\frac{2^{n-1}n \pi^{2n}}{4^n-1}\sum_{(V,P)}\frac{1}{2^{\sharp P}}\prod_{p \in P\setminus\{n\}} { N_{PV}(p) \choose 2} \prod_{r\in [n]\setminus (P\cup V)}N_{PV}(r)\prod_{m=1}^n\frac{1}{2m+N_{PV}(m+1)}.$$

\noindent Let
\begin{equation}\label{4}
\frac{x}{e^x-1}:=\sum_{k=1}^{\infty}\frac{B_n}{n!}x^n=1-\frac{x}{2}+\frac{x^2}{6}-\frac{x^4}{30}+\frac{x^6}{42}-\cdots,
\end{equation}
where $B_n$ are called the Bernoulli numbers. Now applying the Euler formula (cf. \cite[p. 28]{AGHKLT-21})
\[\zeta(2n)=\frac{(-1)^{n+1}(2\pi)^{2n}B_{2n}}{2(2n)!},\]
we obtain the following formula for the Bernoulli number.
\begin{corollary}\label{Bern}
For any integer $n>1$, we have
\begin{equation*}
	B_{2n}=\frac{(-1)^{n+1}}{2^n(2^{2n}-1)}\sum_{(V,P)}\frac{1}{2^{\sharp P}}\prod_{p \in P\setminus\{n\}} { N_{PV}(p) \choose 2} \prod_{r\in [n]\setminus (P\cup V)}N_{PV}(r)\prod_{m=1}^n\frac{1}{2m+N_{PV}(m+1)},
\end{equation*}
where the  pair $(V,P)$ runs over all the $n$-admissible pairs.
\end{corollary}

Bernoulli numbers have many important applications in number theory, such as in the proof of Fermat's Last Theorem. 
The definition of Bernoulli function \eqref{4} is defined via the generating function $\frac{x}{e^x-1}$. 
In general, we use recursion to obtain the representation of Bernoulli numbers.
However, in Corollary \ref{Bern}, we provide a new closed combinatorial formula for Bernoulli numbers
using the widely studied combinatorial numbers $c_n(V.P)$.

We end this section with concrete calculations for $n\leq 4$. As $n=2$, we have $V=\{1\}$, $P=\{2\}$, and $c_2(V,P)=1$, hence
$$
\zeta(4)=\frac{2\pi^{4}}{4^2-1}\frac{1}{2\cdot 1+1}\frac{1}{2\cdot 2}=\frac{\pi^{4}}{90}.
$$
As $n=3$, we have $V=\{1\}$, $P=\{3\}$, and $c_3(V,P)=2$, hence
$$
\zeta(6)=\frac{3\pi^{6}}{4^3-1} \cdot 2\cdot  \frac{1}{2\cdot 1+1}\frac{1}{2\cdot 2+1}\frac{1}{2\cdot 3}=\frac{\pi^{6}}{945}.
$$

For $n=4$, we have two admissible pairs. One is  $V=\{1\}$, $P=\{4\}$,  $c_4(V,P)=4$,  and the other is $V=\{1,2\}$, $P=\{3,4\}$,  $c_4(V,P)=2$.
Evaluating  integrals in each class separately,  we get
\begin{align*}
\zeta(8) &=\frac{4\pi^{8}}{4^4-1} \cdot
\left( 2\cdot  \frac{1}{2\cdot 1+1}\frac{1}{2\cdot 2+2}\frac{1}{2\cdot 3 +1}\frac{1}{2\cdot 4}
+4\cdot  \frac{1}{2\cdot 1+1}\frac{1}{2\cdot 2+1}\frac{1}{2\cdot 3+1}\frac{1}{2\cdot 4}\right)\\
&=\frac{4\pi^{8}}{4^4-1} \cdot \frac{17}{2520}=\frac{\pi^8}{9450}=\frac{2^{8-1}}{8!}|B_8|\pi^8.
\end{align*}

\noindent 
{\bf Acknowledgements.}
The authors gratefully  thank  the anonymous referees for their comments and 
suggestions, which definitely helped to improve the readability and quality of this manuscript.
The authors would like to thank Fritz Gesztesy for pointing the reference \cite{AGHKLT-21}, and  Jun Ma for 
helpful comments and suggestions.

The third named author is supported by the National Natural Science Foundation of China (Grant 12071371).

\end{document}